\documentclass{amsart}

\usepackage{anysize}
\usepackage[latin1]{inputenc}

\marginsize{3cm}{3cm}{3.5cm}{3.5cm}

\newcommand{\jacobi}[2]{\genfrac{(}{)}{.8pt}{}{#1}{#2}}

\begin{document}

\title{Mosaic supercongruences of Ramanujan-type}
\author{Jesús Guillera}
\email{jguillera@gmail.com}
\address{Av. Cesáreo Alierta, 31 esc. izda {\rm $4^o$}--A, Zaragoza (Spain)}
\date{}

\maketitle

\begin{abstract}
We generalize the patterns of supercongruences of Ramanujan-type observed by L. Van Hamme and W. Zudilin to series involving simple square roots anywhere and not only in the result of the sum. To support our observations we give some examples.
\end{abstract}

\section{Ramanujan and Ramanujan-Sato series}

The research of Srinivasa Ramanujan on elliptic integrals and modular equations lead him to discover $17$ surprising series for $1/\pi$ published in $1914$. They are of the following form
\[ \sum_{n=0}^{\infty} A_n (a+bn) z^n = \frac{1}{\pi}, \]
where $z$, $a$ and $b$ are algebraic numbers and
\[
A_n=\frac{\left( \frac{1}{2} \right)_n (s)_n (1-s)_n}{(1)_n^3}, \, \, s=1/2, \, 1/4, \, 1/3, \, 1/6.
\]
During a long time these Ramanujan's series for $1/\pi$ were almost ignored but since $1987$ the interest of mathematicians in them was great. The $17$ formulas as well as many other series of the same type are already proved rigorously \cite{BaBeCh}. Naturally, they are called Ramanujan-type series. In $2002$ T. Sato surprised the mathematical community presenting a series like those of Ramanujan but involving the Apéry numbers
\[ A_n=\sum_{k=0}^n \binom{n}{k}^2\binom{n+k}{k}^2. \]
Inspired in this result mathematicians discovered similar series involving other kinds of special numbers as the Domb numbers
\[ A_n=\sum_{k=0}^n \binom{n}{k}^2\binom{2k}{k}\binom{2n-2k}{n-k}, \]
the Almkvist-Zudilin numbers
\[ A_n=\sum_{k=0}^n (-1)^{n-k} \frac{3^{n-3k} (3k)!}{(k!)^3} \binom{n}{3k}\binom{n+k}{k}, \]
or others, like for example
\[
A_n=\sum_{k=0}^n \binom{2k}{k}^2 \binom{2n-2k}{n-k}^2, \quad A_n=\binom{2n}{n} \sum_{k=0}^n \binom{n}{k}^3,
\quad  A_n=\binom{2n}{n} \sum_{k=0}^n \binom{n}{k}^2 \binom{2k}{k}.
\]
The series involving these kind of numbers are called of Ramanujan-Sato-type and we will refer to these numbers $A_n$ as Ramanujan-Sato-type sequences of numbers. We have taken the definitions from \cite{AlZu}, \cite{ChVe},  \cite{ChTaYaZu}.

\section{Mosaic supercongruences}

We generalize the patterns of supercongruences of Ramanujan-type observed by Van Hamme in \cite{VH} and W. Zudilin in \cite{Zu} to series involving simple square roots anywhere and not only in the result of the sum.

\subsection*{The new pattern}
Suppose $A_n$ is a Ramanujan-Sato sequence of numbers, $z$, $a$, $b$ are algebraic numbers, and
\[ \sum_{n=0}^{\infty} A_n (a+bn) z^n = \frac{1}{\pi}, \]
Suppose also that the partial sums belong to $\mathbb{Q}(\sqrt{d_1}, \dots, \sqrt{d_j})$, where  $d_1,\dots,d_j$ are square free integers, that is
\[ \sum_{n=0}^{p-1} A_n (a+bn) z^n =\alpha_1(p) \sqrt{d_1}+\cdots+\alpha_j(p) \sqrt{d_j}, \]
where $\alpha_1(p),\dots,\alpha_j(p)$ are rational and $a=a_1 \sqrt{d_1}+a_2 \sqrt{d_2}+\cdots+a_j \sqrt{d_j}$.
Then, for primes $p>p_0$, we have the supercongruences
\[ \alpha_i(p) \equiv a_i \jacobi{-d_i}{p} p \pmod{p^3} \qquad i=1,2,\dots, j. \]
We will refer to them as \textit{mosaic supercongruences} because they are as pieces of a single sum expression. We recall that for some kinds of Ramanujan-Sato-type numbers the supercongruences hold only $\! \! \pmod{p^2}$. \par For the Ramanujan-like series for $1/\pi^2$, discovered by the author, we conjecture analogue mosaic supercongruences again generalizing Zudilin's pattern. See the last two examples.

\section{Examples}
All the congruences in the following examples remain unproven.

\subsection*{Example 1}
For the Ramanujan-type series
\[
\frac{\sqrt{15}}{2^7\cdot 5^2} \cdot \sum_{n=0}^{\infty} {\left( \frac{1}{2} \right)_n \left(\frac{1}{6} \right)_n
\left( \frac{5}{6} \right)_n \over (1)_n^3}  (263+5418n) \frac{(-1)^n}{80^{3n}} =\frac{1}{\pi},
\]
we have checked that if we write
\[
\sqrt{15}\sum_{n=0}^{p-1} {\left( \frac{1}{2} \right)_n \left(\frac{1}{6} \right)_n
\left( \frac{5}{6} \right)_n \over (1)_n^3} (263+5418n) \frac{(-1)^n}{80^{3n}}=\alpha_p \sqrt{15},
\]
then, for primes $p>5$, we have the following supercongruences
\[ \alpha_p \equiv 263 \jacobi{-15}{p} p \pmod{p^3}, \]
that is
\[
\sum_{n=0}^{p-1} {\left( \frac{1}{2} \right)_n \left(\frac{1}{6} \right)_n
\left(\frac{5}{6} \right)_n \over (1)_n^3} (263+5418n) \frac{(-1)^n}{80^{3n}} \equiv 263 \jacobi{-15}{p} p \pmod{p^3} \quad p>5,
\]
which is \cite[eq. 21]{Zu}.

\subsection*{Example 2}
For the Ramanujan-type series
\[
\sum_{n=0}^{\infty} \frac{\left( \frac{1}{2} \right)_n \left( \frac{1}{3} \right)_n
\left( \frac{2}{3} \right)_n}{(1)_n^3} \left( \frac{7 \sqrt{7}-10}{27} +  \frac{13 \sqrt{7}-7}{9} n \right) \left( \frac{13 \sqrt 7 -34}{54} \right)^n =\frac{1}{\pi},
\]
we have checked that if we write
\[
\sum_{n=0}^{p-1} \frac{\left( \frac{1}{2} \right)_n \left(\frac{1}{3} \right)_n
\left( \frac{2}{3} \right)_n}{(1)_n^3} \left[\frac{}{} (7 \sqrt{7}-10) + (39 \sqrt{7}-21)n \right]
\left( \frac{13 \sqrt 7 -34}{54} \right)^n =\alpha_p+\beta_p \sqrt{7},
\]
then, for primes $p>7$, we have the following supercongruences
\[ \alpha_p \equiv -10 \jacobi{-1}{p} p \quad \beta_p \equiv 7 \jacobi{-7}{p} p \pmod{p^3}. \]

\subsection*{Example 3}
Consider the Apéry sequence
\[ A_n=\sum_{k=0}^n \binom{n}{k}^2\binom{n+k}{k}^2, \qquad n=0,1,2,\dots\,. \]
Sato's 2002 formula says that
\[
\sum_{n=0}^\infty A_n\left[ \frac{}{} (60\sqrt{15}-134\sqrt{3})+(72\sqrt{15}-160\sqrt{3})n\right]
\biggl(\frac{\sqrt{5}-1}{2}\biggr)^{12n}
=\frac{1}{\pi}.
\]
If we write
\[
\sum_{n=0}^{p-1} A_n\left[ \frac{}{} (60\sqrt{15}-134\sqrt{3})+(72\sqrt{15}-160\sqrt{3})n\right]
\biggl(\frac{\sqrt{5}-1}{2}\biggr)^{12n}
=\alpha_p \sqrt{3}+\beta_p \sqrt{15},
\]
then, for primes $p>5$, we have the supercongruences
\[ \alpha_p \equiv -134 \jacobi{-3}{p}p \quad \beta_p \equiv 60 \jacobi{-15}{p} p \pmod{p^3}. \]

\subsection*{Example 4}
The Ramanujan-type series in \cite[eq 6.1]{BaBe} is
\begin{multline}\nonumber
\sum_{n=0}^{\infty} \frac{\left( \frac{1}{2} \right)_n \left(\frac{1}{6} \right)_n
\left( \frac{5}{6} \right)_n}{(1)_n^3}
\left[\frac{}{}(73+52\sqrt{2}-42\sqrt{3}-30\sqrt{6})+(168+120\sqrt{2}-96\sqrt{3}-69\sqrt{6})n\right] \times \\ (-18872-13344\sqrt{2}+10896\sqrt{3}+7704\sqrt{6})^n=\frac{1}{\pi}. \nonumber
\end{multline}
It derives also from Ramanujan series \cite[eq 6.4]{BaBe} by Zudilin's translation method. Writing
\begin{multline}\nonumber
\sum_{n=0}^{p-1} \frac{\left( \frac{1}{2} \right)_n \left(\frac{1}{6} \right)_n
\left( \frac{5}{6} \right)_n}{(1)_n^3}
\left[\frac{}{}(73+52\sqrt{2}-42\sqrt{3}-30\sqrt{6})+(168+120\sqrt{2}-96\sqrt{3}-69\sqrt{6})n\right] \times \\ (-18872-13344\sqrt{2}+10896\sqrt{3}+7704\sqrt{6})^n= \alpha_p+\beta_p \sqrt{2}+\gamma_p\sqrt{3}+\delta_p\sqrt{6}, \nonumber
\end{multline}
we have, for primes $p>3$, the supercongruences
\[
\alpha_p \equiv 73 \jacobi{-1}{p} p \quad \beta_p \equiv 52 \jacobi{-2}{p} p \quad \gamma_p \equiv -42 \jacobi{-3}{p} p \quad \delta_p \equiv -30 \jacobi{-6}{p} p \pmod{p^3}.
\]

\subsection*{Example 5}
The "complex" Ramanujan series \cite[eq. 45]{GuiZu} is
\[
\sum_{n=0}^{\infty}\frac{(\frac12)_n^3}{(1)_n^3} \biggl(\frac{7\sqrt{7}-13\sqrt{-1}}{64}+\frac{15\sqrt{7}-21 \sqrt{-1}}{32}n \biggr) \biggl(\frac{47+45\sqrt{-7}}{128}\biggr)^n=\frac{1}{\pi}.
\]
Writing
\[
\sum_{n=0}^{p-1}\frac{(\frac12)_n^3}{(1)_n^3} \left[ \frac{}{} (7\sqrt{7}-13\sqrt{-1})+(30\sqrt{7}-42 \sqrt{-1})n \right] \biggl(\frac{47+45\sqrt{-7}}{128}\biggr)^n=\alpha_p\sqrt{-1}+\beta_p\sqrt{7},
\]
we have, for primes $p>7$, the supercongruences
\[ \alpha_p \equiv -13 p \quad \beta_p \equiv 7 \jacobi{-7}{p} p \pmod{p^3}.\]

\subsection*{Example 6}
For the Ramanujan-like series \cite[eq. 10]{coleccion}
\[
\sum_{n=0}^{\infty}  \frac{\left(\frac12 \right)_n \left( \frac13 \right)_n \left( \frac23 \right)_n \left( \frac16 \right)_n \left( \frac56 \right)_n}{(1)_n^5} (-1)^n \left( \frac{3}{4} \right)^{6n} (45+549n+1930n^2)=\frac{384}{\pi^2},
\]
we have checked that if we write
\[
\sum_{n=0}^{p-1}  \frac{\left(\frac12 \right)_n \left( \frac13 \right)_n \left( \frac23 \right)_n \left( \frac16 \right)_n \left( \frac56 \right)_n}{(1)_n^5} (-1)^n \left( \frac{3}{4} \right)^{6n} (45+549n+1930n^2)=\alpha_p,
\]
then, for primes $p>3$, we have the supercongruences
\[ \alpha_p \equiv 45 p^2 \pmod{p^5}, \]
that is, they follow Zudilin's pattern \cite{Zu0}.

\subsection*{Example 7}
The only known (unproven) hypergeometric Ramanujan-like series for $1/\pi^2$ with a non rational value of $z$ is
\begin{multline}\nonumber
\sum_{n=0}^{\infty} \frac{\left(\frac12 \right)_n^3 \left( \frac13 \right)_n \left( \frac23 \right)_n}{(1)_n^5} \left( \frac{15\sqrt{5}-33}{2} \right)^{3n} \times \\ \nonumber
\left[ \frac{}{} (56-25\sqrt{5})+(303-135\sqrt{5})n+(1220/3-180\sqrt{5})n^2\right]=\frac{1}{\pi^2}, \nonumber
\end{multline}
see \cite[eq. 9]{coleccion}. Writing
\begin{multline}\nonumber
\sum_{n=0}^{p-1}  \frac{\left(\frac12 \right)_n^3 \left( \frac13 \right)_n \left( \frac23 \right)_n}{(1)_n^5} \left( \frac{15\sqrt{5}-33}{2} \right)^{3n} \times \\ \nonumber
\left[ \frac{}{} (56-25\sqrt{5})+(303-135\sqrt{5})n+(1220/3-180\sqrt{5})n^2\right]=\alpha_p+\beta_p\sqrt{5}, \nonumber
\end{multline}
we have, for primes $p>5$, the supercongruences
\[ \alpha_p \equiv 56 p^2 \quad \beta_p \equiv -25 \jacobi{5}{p} p^2 \pmod{p^5}, \]
which generalizes Zudilin's pattern.

\subsection*{More support}
To support even more our observations, we have considered many of the series, involving only simple square roots in \cite{BaBe}, \cite{BoBo} and \cite{ChVe}, and observed the expected mosaic supercongruences in all the cases.

\section{Concluding remarks}

An excellent survey on Ramanujan-type series is \cite{BaBeCh} and a beautiful survey on recent advances is in \cite{Zu0}. There are many examples of convergent Ramanujan-type and Ramanujan-Sato-type series in the literature, the most spectacular are of simple and very fast series. From the modular theory of Ramanujan-type series we know that there are functions $z(q)$, $b(q)$ and $a(q)$, with $q=e^{i \pi \tau}$ and ${\rm Im}(\tau)>0$, which take algebraic values when $\tau$ is a quadratic irrational. Obviously the series are faster as ${\rm Im}(\tau)$ is bigger. If ${\rm Im}(\tau)$ is small  enough then the series is "divergent". An example of a "divergent" Ramanujan-type series is in \cite[p. 371]{BoBo}, which corresponds to ${\rm Im}(\tau)=\sqrt{253}/11$. Convergent or divergent series lead to supercongruences following exactly the same patterns \cite{GuiZu}. This fact is not very surprising if we observe that the sum from $n=0$ to $p-1$ do not discern the convergent or divergent origin.
\par Taking into account that the Jacobi symbols are the quadratic residues, perhaps it can give the clue to discover the pattern when the algebraic numbers involved are more complicated. We will continue investigating on this idea.

\subsection*{Acknowledgment}
I thank to Wadim Zudilin for encourage me to discover a pattern when the series involve algebraic numbers.


\begin{thebibliography}{99}

\bibitem{AlZu}
\textsc{G. Almkvist} and \textsc{W. Zudilin},
Differential equations, mirror maps and zeta values,
in Mirror Symmetry V, N. Yui, S.-T. Yau, and J.D. Lewis (eds.), AMS/IP Studies in Advanced Mathematics 38 (2007), International Press \& Amer. Math. Soc., 481--515
(e-print math.NT/0402386)

\bibitem{BaBeCh}
\textsc{N.D. Baruah}, \textsc{B.C. Berndt}, \textsc{H.H. Chan},
Ramanujan's series for $1/\pi$: A survey,
\emph{The Amer. Math. Monthly} \textbf{116}, (2009), 567--587;
also available at the pages http://www.math.ilstu.edu/cve/speakers/Berndt-CVE-Talk.pdf, and http://www.math.uiuc.edu/~berndt/articles/monthly567-587.pdf.

\bibitem{BaBe}
\textsc{N.~Baruah} and \textsc{B.~Berndt},
Eisenstein series and Ramanujan-type series for $1/\pi$,
\emph{Preprint} (2009); also available at
http://citeseerx.ist.psu.edu/viewdoc/summary?doi=10.1.1.157.8033

\bibitem{BoBo}
\textsc{J.~Borwein} and \textsc{P.~Borwein},
More Ramanujan-type series for $1/\pi$, in Ramanujan Revisited, G. E. Andrews et al. (eds.) (Academic Press, Boston, 1988), 359--374; \\ also available at http://www.cecm.sfu.ca/personal/pborwein/PAPERS/CP4.pdf

\bibitem{ChVe}
\textsc{Heng Huat Chan} and \textsc{H.~Verril},
The Apéry numbers, the Almkvist--Zudilin numbers and new series for $1/\pi$.
\emph{Math. Res. Lett.} \textbf{16}:3 (2009), 405--420; also available at
http://www.mathjournals.org/mrl/2009-016-003/2009-016-003-003.pdf

\bibitem{ChTaYaZu}
\textsc{Heng Huat Chan}, \textsc{Y. Taniyawa}, \textsc{Y. Yang} and \textsc{W. Zudilin},
New analogues of Clause's identities arising from the theory of modular forms

\bibitem{GuiZu}
\textsc{J.~Guillera} and \textsc{W.~Zudilin},
"Divergent" Ramanujan-type supercongruences.
(e-print arXiv:1004.4337)

\bibitem{coleccion}
\textsc{J.~Guillera},
Collection of Ramanujan-like series for $1/\pi^2$. Unpublished manuscript available at the author's web site.

\bibitem{VH}
\textsc{L. Van Hamme},
Some conjectures concerning partial sums of generalized hypergeometric
series, \textit{p-adic functional analysis} (Nijmegen, 1996), Lecture Notes in Pure and Appl. Math.
192, Dekker, New York (1997), 223--236.

\bibitem{Zu0}
\textsc{W.~Zudilin},
Ramanujan-type formulae for $1/\pi$: A second wind?
(Banff, June 3--8, 2006), N. Yui, H. Verrill, \& C.F. Doran (eds.), Fields Inst. Commun. Ser. 54 (2008), Amer. Math. Soc. and Fields Inst., 179--188 (e-print arXiv:0712.1332 [math.NT])

\bibitem{Zu}
\textsc{W.~Zudilin},
Ramanujan-type supercongruences,
\emph{J. Number Theory} \textbf{129}:8 (2009), 1848--1857 (e-print arXiv:0805.2788 [math.NT])

\end{thebibliography}
\end{document}